\documentstyle[12pt]{article}

\font\teneufm=eufm10 \font\seveneufm=eufm7 \font\fiveeufm=eufm5
\newfam\eufmfam
\textfont\eufmfam=\teneufm \scriptfont\eufmfam=\seveneufm
\scriptscriptfont\eufmfam=\fiveeufm
\def\frak#1{{\fam\eufmfam\relax#1}}

\newfam\msbfam
\font\tenmsb=msbm10 scaled \magstep1   \textfont\msbfam=\tenmsb
\font\sevenmsb=msbm7 scaled \magstep1 \scriptfont\msbfam=\sevenmsb
\font\fivemsb=msbm5 scaled \magstep1
\scriptscriptfont\msbfam=\fivemsb
\def\Bbb{\fam\msbfam \tenmsb}

\def\RR{{\Bbb R}}
\def\CC{{\Bbb C}}
\def\QQ{{\Bbb Q}}
\def\NN{{\Bbb N}}
\def\ZZ{{\Bbb Z}}
\def\II{{\Bbb I}}
\def\11{{\bf 1}}

\def\e{\epsilon}
\def\q{{\bf q}}
\def\v{{\bf v}}
\def\ra{\rightarrow}
\def\ss{\subseteq}
\def\e{\epsilon}
\def\F{{\cal F}}
\def\O{\Omega}
\def\bO{\partial \Omega}
\def\A{{\cal A}}
\def\endpf{\hfill $\Box$ \smallskip \\}
\def\a{\alpha}
\def\tg{\widetilde{\gamma}}
\def\sm{\setminus}
\def\ss{\subseteq}

\newtheorem{definition}{\bf Definition}[section]
\newtheorem{theorem}[definition]{\bf Theorem}

\newtheorem{remark}[definition]{\sc Remark}
\newtheorem{proposition}[definition]{\bf Proposition}
\newtheorem{example}[definition]{\sc Example}

\begin{document}

\begin{center}
\Large \bf Quantum Normal Families:
\smallskip \\
\large Normal Families of Holomorphic \\
Functions and Mappings\\
on a Banach Space\footnote{Both authors were guests at the
American Institute of Mathematics (AIM) during a portion of this
work.  AIM sponsored a workshop on holomorphic mappings that was
particularly useful to these studies.}
\end{center}
\begin{center}
by \\
Kang-Tae Kim,\footnote{Author supported in part by grant
R01-1999-00005 from The Korean Science and Engineering
Foundation.} and Steven G. Krantz,\footnote{Author supported in
part by grant DMS-9988854 from the National Science Foundation.}
\end{center}

\begin{quotation}
\noindent \small \bf Abstract:  \sl The authors lay the
foundations for the study of normal families of holomorphic
functions and mappings on an infinite-dimensional normed linear
space.  Characterizations of normal families, in terms of value
distribution, spherical derivatives, and other geometric
properties are derived.  Montel-type theorems are established.

A number of different topologies on spaces of holomorphic mappings
are considered.  Theorems about normal families are formulated and
proved in the language of these various topologies.

Normal functions are also introduced.  Characterizations in terms
of automorphisms and also in terms of invariant derivatives are
presented.
\end{quotation}

\setcounter{section}{-1}

\section{Introduction}

The theory of holomorphic functions of infinitely many complex
variables is about forty years old.  Pioneers of the subject were
Nachbin [NAC1--NAC17], Gruman and Kiselman [GRK], and Mujica
[MUJ]. After a quiet period of nearly two decades, the discipline
is now enjoying a rebirth thanks to the ideas of Lempert
[LEM1--LEM6]. Lempert has taught us that it is worthwhile to
restrict attention to {\it particular} Banach spaces, and he has
directed our efforts to especially fruitful questions.

The work of the present paper is inspired by the results of [KIK].
That paper studied domains in a Hilbert space with an automorphism
group orbit accumulating at a boundary point.  As was the case in
even one complex variable, normal families played a decisive role
in that study.  With a view to extending those explorations, it
now seems appropriate to lay the foundations for normal families
in infinitely many complex variables.

One of the thrusts of the present paper is to demonstrate that
normal families may be understood from several different points of
view.  These include:
\begin{enumerate}
\item Classical function theory
\item Hyperbolic geometry
\item Functional analysis
\item Distribution theory
\item Currents
\item Comparison of different topologies and norms on the space of
holomorphic functions
\end{enumerate}
It is our intention to explain these different approaches to the
subject and to establish relationships among them.

A second thrust is to relate the normality of a family on the
entire space $X$ (or on a domain in $X$) to the normality of the
restriction of the family to slices (suitably formulated).  This
point of view was initiated in [CIK], and it has proved useful and
intuitively natural.

Throughout this paper, $X$ is a separable Banach space, $\O$ is a
domain (a connected open set) in $X$, $\F = \{f_\a\}_{\a \in {\cal
A}}$ is a family of holomorphic functions on $\O$, and $D \ss \CC$
is the unit disc.  If $\O'$ is another domain in some other
separable Banach space $Y$, then we will also consider families
$\{F_\alpha\}$ of holomorphic mappings from $\O$ to $\O'$.
Although separability of $X$ is not essential to all of our
results, it is a convenient tool in many arguments.  Certainly, in
the past, the theory of infinite dimensional holomorphy has been
hampered by a tendency to shy away from such useful extra
hypotheses.

Part of the beauty and utility of studying infinite-dimensional
holomorphy is that the work enhances our study of
finite-dimensional holomorphy.  Indeed, it is safe to say that the
present study has caused us to re-invent what a normal family of
holomorphic functions ought to be.

One of the interesting features of the present work, making it
different from more classical treatments in finite dimensions, is
that compact sets now play a different role.  If $W$ is a given
open set in our space $X$ (say the unit ball in a separable
Hilbert space), then $W$ cannot be exhausted by an increasing
union of compact sets in any obvious way.  Another feature is
that, in finite dimensions, all reasonable topologies on the space
of holomorphic functions on a given domain are equivalent.  In
infinitely many variables this is no longer the case, and we hope
to elucidate the matter both with examples and results relating
the different topologies.

It is a pleasure to thank John M\raise3pt\hbox{c}Carthy for
helpful conversations about various topics in this paper.

\section{Basic Definitions}

We will now define holomorphic functions and mappings and normal
families.  We refer the reader to the paper [KIK] and the book
[MUJ] for background on complex analysis in infinite dimensions.

\begin{definition} \rm
A {\it domain} $\Omega \subseteq X$ is a connected open set.
\end{definition}

\let\foo\it
\begin{definition}  \rm
Let $\Omega \subseteq X$ be an open set.  Let $u: \Omega
\rightarrow Y$ be a mapping, where $Y$ is some other separable
Banach space. For ${\bf q} \in \Omega$ and $\v_1,\dots, \v_k \in
X$, we define the derivatives
$$
du({\bf q}; \v_j) = \lim_{\RR \ni \epsilon \ra 0} \frac{u(\q +
\epsilon \v_j) - u(\q)}{\e} \, .
$$
and
$$
\overline{D} u (\q; \v) = \frac{du(\q;\v) + i du(\q; i\v)}{2} \, .
$$
\end{definition}
A function $f$ on $\Omega$ is said to be {\it continuously
differentiable}, or $C^1$, if $df({\bf q;\v})$ exists for every
point ${\bf q} \in \Omega$ and every vector ${\bf v}$, and if the
resulting function $({\bf q,\v}) \mapsto df({\bf q;\v})$ is
continuous.

\begin{definition} \rm
Let $\Omega \ss X$ be an open set and $f$ a $C^1$-smooth function
or mapping defined on $\Omega$.  We say that $f$ is {\it
holomorphic} on $\Omega$ if $\overline{D} f \equiv 0$ on $\Omega$.
\end{definition}

The definition just given of ``holomorphic function'' or
``holomorphic mapping'' is equivalent, in the $C^1$ category, to
requiring that the restriction of the function or mapping to every
complex line be holomorphic in the classical sense of the function
theory of one complex variable. We shall have no occasion, in the
present paper, to consider functions that are less than $C^1$
smooth, but holomorphicity can, in principle, be defined for
rougher functions.

\let\foo\it
\begin{definition} \rm
Let ${\cal F} = \{f_\alpha\}_{\alpha \in {\cal A}}$ be a family of
holomorphic functions on a domain $\Omega \ss X$. We say that
${\cal F}$ is a {\foo normal family} if every subsequence $\{f_j\}
\ss \F$ either
\begin{description}
\item[{\bf 1.4.1 (normal convergence)}]  has a subsequence that
converges uniformly on compact subsets of $\Omega$;

\noindent \hspace*{-.45in} or

\item[{\bf 1.4.2 (compact divergence)}]  has a subsequence $\{f_{j_k}\}$
such that, for each compact $K \ss \Omega$ and each compact $L \ss
\CC$, there is a number $N$ so large that if $k \geq N$ then
$f_{j_k}(K) \cap L = \emptyset$.
\end{description}
\end{definition}

It is convenient at this juncture to define a type of topology
that will be of particular interest for us.  If $\Omega \ss X$ is
a domain and ${\cal O}(\Omega)$ the space of holomorphic functions
on $\Omega$, then we let ${\cal B}$ denote the topology on ${\cal
O}(\Omega)$ of uniform convergence on compact sets. Of course a
sub-basis for the topology ${\cal B}$ is given by the sets (with
$\epsilon > 0$, $g \in {\cal O}(\Omega)$ arbitrary, and $K \ss
{\cal O}$ a compact set)
$$
B_{g,K,\epsilon} = \{ f \in {\cal O}(\Omega) :
                \sup_{z \in K} |f(z) - g(z)| < \epsilon\} \, .
$$
It is elementary to verify (or see [MUJ]) that the limit of a
sequence of holomorphic functions on $\Omega$, taken in the
topology ${\cal B}$, will be another holomorphic function on
$\Omega$.

\let\foo\it
\begin{definition} \rm
Let $\Omega \ss X$ be a domain.  Let ${\cal U} \equiv
\{U_\alpha\}_{\alpha \in {\cal A}}$ be a semi-norm topology on the
space ${\cal O}(\Omega)$ of holomorphic functions on $\Omega$.  We
say that ${\cal U}$ is a {\foo Montel topology} on ${\cal
O}(\Omega)$ if the mapping
\begin{eqnarray*}
\hbox{\rm id}: [{\cal O}(\Omega), {\cal U}] & \longrightarrow &
[{\cal O}(\Omega), {\cal B}] \\
f & \longmapsto     & f
\end{eqnarray*}
is a compact operator.
\end{definition}

\begin{example} \rm
\begin{enumerate}
\item[\hbox{\bf (1)}]  Let $\Omega \ss \CC^n$ be a domain in
finite-dimensional complex space.  The topology ${\cal B}$ is a
Montel topology. This is the content of the classical Montel
theorem on normal families (see [MUJ]).
\item[\hbox{\bf (2)}]  Let $\Omega \ss \CC$ be a domain in one-dimensional
complex space.  Let $k$ be a positive integer, and let superscript
$(k)$ denote the $k^{\rm th}$ derivative. The topology ${\cal
C}_k$ with sub-basis given by the union of the sets (with
$\epsilon > 0$, $g \in {\cal O}(\Omega)$ arbitrary, and $K \ss
{\cal O}$ a compact set)
$$
N_{g,K,\epsilon} = \{ f \in {\cal O}(\Omega) :
                \sup_{z \in K} |f^{(j)}(z) - g^{(j)}(z)| < \epsilon\} \, ,
$$
for $j = 0, 1, \dots, k$, is a Montel topology.  Of course, by
integration (and using the Cauchy estimates), the topology ${\cal
C}_k$ is equivalent to the topology ${\cal B}$.  [A similar
topology can be defined for holomorphic functions on a domain in
the finite-dimensional space $\CC^n$.]
\item[\hbox{\bf (3)}]  Let $\Omega \ss \CC$ be a domain in one-dimensional
complex space.  The topology ${\cal D}$ with sub-basis given by
the sets (with $\epsilon > 0$, $g \in {\cal O}(\Omega)$ arbitrary,
and $\widetilde{\gamma} \ss {\cal O}$ the compact image of a
closed curve $\gamma: [0,1] \rightarrow \Omega$)
$$
M_{g,\widetilde{\gamma},\epsilon} = \{ f \in {\cal O}(\Omega) :
                \sup_{z \in \widetilde{\gamma}} |f(z) - g| < \epsilon\}
$$
is a Montel topology, as the reader may verify by using the Cauchy
estimates.  Of course, once again, the maximum principle may be
used to check that the topology ${\cal D}$ is equivalent to the
topology ${\cal B}$.
\item[\hbox{\bf (4)}]  In the reference [NAC18], Leopoldo Nachbin
defined the concept of a seminorm that is ``ported'' by a compact
set.  We review the notion here.  Let $X$ and $Y$ be separable
Banach spaces as usual.  Let $\O \ss X$ be a domain, and let $K
\ss \O$ be a fixed compact subset.  We consider the family ${\cal
H}(\O, Y)$ of holomorphic mappings from $\O$ to $Y$.  A seminorm
$\rho$ on ${\cal H}(\O, Y)$ is said to be {\it ported} by the set
$K$ if, given any open set $V$ with $K \ss V \ss \O$, we can find
a real number $c(V) > 0$ such that the inequality
$$
\rho(f) \leq c(V) \cdot \sup_{x \in V} \| f(x) \|   \eqno (*)
$$
holds for every $f \in {\cal H}(\O, Y)$.
\smallskip \\
\null \ \ \ \ We note that the holomorphic mapping $f$ here need
not be bounded on $V$. What is true, however (and we have noted
this fact elsewhere in the present paper), is that once $\Omega$
and $K$ are fixed then there will exists some open set $V$ as
above on which $f$ is bounded.  So that, for this choice of $V$,
the inequality $(*)$ will be non-trivial.
\smallskip  \\
\null \ \ \ \ Now we use the notion of ``seminorm ported by $K$''
to define a topology on ${\cal H}(\O, Y)$ as follows:  we consider
the topology induced by all seminorms that are ported by compact
subsets of $\O$.  It is to be noted that, in finite dimensions,
this new topology is no different from the standard compact-open
topology.  But in infinite dimensions it is quite different.  As
an example, let $X = Y = \ell_2$, which is of course a separable
Hilbert space. Let a typical element of $\ell_2$ be denoted by
$\{a_j\}_{j=1}^\infty$, and let the $j^{\rm th}$ coordinate be
$z_j$. Let $\O \ss X$ be a domain and let $K \ss \O$ be a compact
set. Consider holomorphic functions $f:\Omega \ra \CC$.  Define a
semi-norm by
$$
\rho^*(f) \equiv \sum_{j=1}^\infty \sup_K \left |
   \frac{\partial f}{\partial z_j} \right | \, .
$$
Then it is clear, by the Cauchy estimates, that $\rho^*$ is
ported.

\ \ \ \ But it is also clear that a typical open set defined by
$\rho^*$ will not contain any non-trivial open set from the
compact-open topology.  Thus this topology is {\it not} Montel. Of
course it is now a simple matter to generate many other
interesting examples of ported seminorms.
\end{enumerate}
\end{example}

Now we have

\begin{theorem} \sl
Let $\F = \{f_\alpha\}_{\alpha \in {\cal A}}$ be a family of
holomorphic functions on a domain $\O \subseteq X$.  Assume that
there is a finite constant $M$ such that $|f_\alpha(x)| \leq M$
for all $f_\alpha \in \F$ and all $z \in \O$. Let $K$ be a compact
subset of $\O$.  Then every sequence in $\F$ has itself a
subsequence that converges uniformly on $K$.
\end{theorem}
\noindent {\bf Proof:}  Of course the hypothesis of uniform
boundedness precludes compact divergence.  So we will verify {\bf
1.4.1}.  Fix a compact subset $K \ss \O$. Then there is a number
$\eta > 0$ such that if $k \in K$ then $B(k, 3\eta) \ss \Omega$.
Select $f_\alpha \in \F$. Now if $k \in K$ and $\ell$ is any point
such that $\|k - \ell\| < \eta$ then we may apply the Cauchy
estimates (on $B(k, 2\eta)$) to the restriction of $f_\alpha$ to
the complex line through $k$ and $\ell$.  We find that the
$f_\alpha$ have bounded directional derivatives. Therefore they
are (uniformly) Lipschitz and form an equicontinuous family of
functions.

As a result of these considerations, the Arzela-Ascoli theorem
applies to the family $\F$ restricted to $K$.  Thus any sequence
in $\F$ has a subsequence convergent on $K$.
\endpf

In practice, it is useful to have a version of Theorem 1.7 that
hypothesizes only uniform boundedness on compact sets.  This is a
tricky point in the infinite-dimensional setting for the following
reason:  Classically (in finite dimensions), one derives this new
result from (the analog of) Theorem 1.7 by taking a compact set $K
\ss \O$ and fattening it up to a slightly larger compact $L \ss
\O$.  Since the family $\F$ is uniformly bounded on $L$, an
analysis similar to the proof of 1.7 may now be performed.  In the
infinite-dimensional setting this attack cannot work, since there
is no notion of fattening up a compact set to a larger compact
set.

Nonetheless, we have several different ways to prove a more
general, and more useful, version of Montel's theorem.  The
statement is as follows.

\begin{theorem}[Montel]  \sl
Let $\F = \{f_\alpha\}_{\alpha \in {\cal A}}$ be a family of
holomorphic functions on a domain $\O$.  Assume that $\F$ is
uniformly bounded on compact sets, in the sense that for each
compact $L \ss \O$ there is a constant $M_L > 0$ such that
$|f_\alpha(z)| \leq M_L$ for every $z \in L$ and every $f_\alpha
\in \F$. Then every sequence in $\F$ has itself a subsequence that
converges uniformly on each compact set $K \ss \O$.  [Note that we
are saying that there is a single sequence that works for every
set $K$.]
Thus ${\cal F}$ is a normal family.
\end{theorem}

\begin{remark} \rm We may rephrase Montel's theorem by saying that
the topology ${\cal B}$ is a Montel topology.
\end{remark}

\noindent {\bf Proof of the Theorem:} Fix a compact set $K \ss
\O$. Of course the family ${\cal F}$ is bounded on $K$ by
hypothesis. We claim that ${\cal F}$ is bounded on some
neighborhood $U$ of $K$.  To this end, and seeking a
contradiction, we suppose instead that for each integer $N > 0$
there is a point $x_N \in \O$ such that $\hbox{dist}(x_N, K) <
1/N$ and $|f_\alpha(x_N)| > N$.  Then the set
$$
L = K \cup \{x_N\}_{N=1}^\infty
$$
is compact.  So the family $\F$ is bounded on $L$. But that
contradicts the choice of the $x_N$.

We conclude that, for some $N$, $x_N$ does not exist.  That means
that there is a number $N_0 > 0$ such that the family $\F$ is
uniformly bounded on $U \equiv \{x \in \Omega: \hbox{dist}(x, K) <
1/N_0\}$. As a result, we may imitate the proof of Theorem 1.7,
merely substituting $U$ for $\O$. \hfill $\Box$

\begin{remark} \rm
We thank Laszlo Lempert for the idea of the proof of 1.8 just
presented.
\end{remark}

We now indulge in a slight digression, partly for interest's sake
and partly because the argument will prove useful below.  In fact
we will provide a proof of Theorem 1.8 that depends on the
Banach-Alaoglu theorem.  This is philosophically appropriate, for
it validates in yet another way that a normal families theorem is
nothing other than a compactness theorem.  After that we will
sketch a proof that depends on the theory of currents.
\smallskip \\

\noindent {\bf Alternative (Banach-Alaoglu) Proof of Theorem 1.8}:  \\
For clarity and simplicity, we begin by presenting this proof in
the complex plane $\CC$.  The reader who has come this far will
have no trouble adapting the argument to {\it finitely many}
complex variable space $\CC^n$.  We provide a separate argument
below for the infinite dimensional case.

Now fix a domain $\O \ss \CC$.  Let $\F = \{f_\alpha\}_{\alpha \in
{\cal A}}$ be a family of holomorphic functions on $\O$ which is
bounded on compact sets.  Fix a piecewise $C^1$ closed curve
$\gamma: [0,1] \ra \O$.  Let $\widetilde{\gamma}$ denote the {\it
image} of $\gamma$, which is of course a compact set in $\Omega$.
Consider the functions
$$
g_\alpha \equiv f_\alpha \bigr |_{\widetilde{\gamma}} .
$$
Then each $g_\alpha$ is smooth on $\widetilde{\gamma}$ and the
family ${\cal G} \equiv \{g_\alpha\}_{\alpha \in {\cal A}}$ is
bounded by some constant $M$.  So we may think of ${\cal G} \ss
L^\infty(\widetilde{\gamma})$ as a bounded set.  Since
$L^\infty(\tg)$ is the dual of $L^1(\tg)$, we may apply the
Banach-Alaoglu theorem to extract a subsequence (which we denote
by $\{g_j\}$ for convenience) that converges in the weak-$*$
topology.  Call the weak-$*$ limit function $g$.

Now fix a point $z$ that lies in the interior, bounded component
of the complement of $\tg$.  Of course the function
$$
t \longmapsto \frac{\gamma'(t)}{\gamma(t) - z}
$$
lies in $L^1(\tg)$.  So, by weak-$*$ convergence and the Cauchy
integral formula, we know that
$$
g_j(z) \equiv \frac{1}{2\pi i} \oint_\gamma
\frac{g_j(\zeta)}{\zeta - z} \, d\zeta  \rightarrow \frac{1}{2\pi
i} \oint_\gamma \frac{g(\zeta)}{\zeta - z} \, d\zeta
    \equiv G(z) \, .
$$
Here the last equality {\it defines} the function $G$.

So we see that the functions $g_j$, which of course must agree
with $f_j$ at points inside the curve $\gamma$, tend pointwise to
the function $G$; and the function $G$ is perforce holomorphic
inside the image curve $\tg$. We will show that in fact the
convergence is uniform on compact sets inside of $\tg$.

So fix a compact set $K$ that lies in the bounded open set
interior to $\tg$.  Fix a piecewise $C^1$, simple, closed curve
$\gamma^*$ whose image is disjoint from, and lies inside of,
$\tg$, and which surrounds $K$.  Let $\eta > 0$ be the distance of
$K$ to $\widetilde{\gamma^*}$, the image of $\gamma^*$. Now fix a
small $\epsilon > 0$ (here $\epsilon$ should be smaller than the
length of $\gamma^*$). Choose a set $E \ss \widetilde{\gamma^*}$
such that $E$ has linear measure less than $\epsilon$ and so that
(by Lusin's theorem)
$$
|g_j(\zeta) - g(\zeta)| < \epsilon
$$
when $j$ is sufficiently large ($j > N$, let us say) and $\zeta
\in \widetilde{\gamma^*} \setminus E$.

Then, for $j, k > N$ and $z \in K$ we have
\begin{eqnarray*}
|g_j(z) - g_k(z)| & \leq & \frac{1}{2\pi} \int_{{}^c E} \left |
   \frac{g_j(\zeta) - g_k(\zeta)}{\zeta - z} \right | d|\zeta|
+ \frac{1}{2\pi} \int_E \left |
   \frac{g_j(\zeta) - g_k(\zeta)}{\zeta - z} \right | d|\zeta| \\
& \leq & \frac{1}{2\pi} \hbox{\rm length}(\widetilde{\gamma^*})
\frac{\epsilon}{\eta}
 + \frac{1}{2\pi} \cdot \epsilon \cdot \frac{2M}{\eta} \, .
\end{eqnarray*}
Since $\epsilon > 0$ may be chosen to be arbitrarily small, we
conclude that $g_j \ra g$ uniformly on the compact set $K$.  That
is what we wished to prove for the single compact set $K$.

We note that this proof may be performed when $\gamma$ is a
positively oriented curve describing {\it any} square inside
$\Omega$ with sides parallel to the axes, rational center, and
rational side length.  Of course it is always possible to produce
the curve $\gamma^*$ as the union of finitely many such curves. As
a result, the usual diagonalization procedure may be formed over
these countably many curves, producing a single subsequence that
converges uniformly on any compact set in $\Omega$ to a limit
function $G$.
\endpf

\noindent {\bf Alternative (Currents) Proof of Theorem 1.8
for Separable Banach Spaces}:  \\
We refer to the very interesting paper [ALM] of Almgren.  That
paper gives a characterization of the dual of the space of all
$k$-dimensional, real, rectifiable currents in $\RR^N$.
Remarkably, Almgren's proof uses both the Continuum Hypothesis and
the Axiom of Choice.  An examination of Almgren's proof reveals
that the arguments are also valid when $\RR^N$ is replaced by any
separable Banach space. We take that result for granted, and leave
it to the reader to check the details in [ALM].

Accepting that assertion, we see that the hypothesis of uniform
boundedness of a family $\F$ of holomorphic functions on compact
subsets of a domain $\O$ in a separable Banach space $X$ can be
interpreted as a boundedness statement about one-dimensional
holomorphic currents.  Specifically, let $\F$ be a family of
holomorphic functions on a domain $\O \ss X$, and assume that $\F$
is bounded on compact subsets of $\O$.  As we have seen (proof of
Theorem 1.8), it follows that if $K \ss \O$ is any compact set
then there is a small neighborhood $U$ of $K$, with $K \subset \!
\subset U \ss \O$, such that $\F$ is bounded on $U$.  As a result,
we may apply Cauchy estimates to see that if $\F =
\{f_\alpha\}_{\alpha \in {\cal  A}}$ then $\F' = \{\partial
f_\alpha\}_{\alpha \in {\cal  A}}$ is bounded on $K$.  But then,
by the generalization of Almgren's theorem to infinite dimensions,
we may think of $\F'$ as a bounded family in the dual of the space
of 1-dimensional (complex) currents on $\O$. By the Banach-Alaoglu
theorem, we may therefore extract from any sequence in $\F'$ a
weak-$*$ convergent subsequence.  Call it, for convenience,
$\{f_j\}$.

But now it is possible to imitate the first alternative proof of
Theorem 1.8 as follows.  Fix a closed, piecewise $C^1$ curve
$\gamma: [0,1] \ra \Omega$ {\it that bounds an analytic disc {\bf
d}} in $\Omega$.  Think of the elements $\partial f_j$ restricted
to the image $\widetilde{\gamma}$ of this curve.  They form a
bounded family in $L^\infty(\widetilde{\gamma})$.  Thus the first
alternative proof may be imitated, step by step, to produce a
limit holomorphic function on the analytic disc {\bf d}. In fact
we may even take the argument a step further.  We may look at any
$k$-dimensional slice of $\Omega$ and use the Bochner-Martinelli
kernel instead of the 1-dimensional Cauchy kernel to find that
there is a uniform limit on any compact subset of any
$k$-dimensional slice of $\Omega$.  This produces the required
limit function $G$ for the subsequence $f_j$.  [Note that, because
we are assuming the space to be separable, we can go further an
even extract a subsequence that converges on {\it every} compact
subset.  More will be said about this point in the next remark and
in what follows.]
\endpf

\begin{remark}  \rm
This last is still not the optimal version of what we usually call
Montel's theorem.  In the classical, finite-dimensional
formulation of Montel's result we usually derive a single
subsequence that converges uniformly on {\it every} compact set.
The question of whether such a result is true in infinite
dimensions is complicated by the observation that it is no longer
possible, in general, to produce a sequence of sets $K_1 \subset
\subset K_2 \subset \subset \cdots \subset \subset X$ for our
Banach space $X$ with the property that each compact subset of $X$
lies in some $K_j$.  In fact the full-bore version of Montel's
theorem, as just described, is false.  The next example of Y. Choi
illustrates what can go wrong.
\end{remark}

\begin{example} \rm
Consider the Banach space $X = \ell^\infty$. Let
$$
e_j = (0, \ldots, 0, 1, 0, \ldots) \, ,
$$
in which all components except the $j^{\rm th}$ are zero.  Let
$e_j^* : X \to \CC$ be defined by
$$
e_j^* \left (\sum_{k=1}^\infty a_k e_k \right ) = a_j.
$$
This function is obviously holomorphic.  However, the sequence
$\{e_j^*\}$ does not have a subsequence that converges uniformly
on compact subsets.  To see this, let us assume to the contrary
that $\{e_{j_m}^*\}_{m=1}^\infty$ is a subsequence that converges
uniformly on compact subsets.  Then in particular it should
converge on singleton set consisting of the point ${\bf p}$ that
is given by
$$
{\bf p} = \sum_{m=1}^\infty (-1)^m e_{j_m} \, .     \eqno{(1)}
$$
But, $e_{j_m}^* ({\bf p}) = (-1)^m$, and this sequence of scalars
does not converge.
\end{example}

It should be noted that this example can be avoided if we demand
in advance that the Banach space $X$ be separable. One simply
produces a countable, dense family of open balls, extracts a
convergent sequence for each such ball, and then diagonalizes as
usual.  Mujica [MUJ], in his treatment of normal families,
achieves the full result by adding a hypothesis of pointwise
convergence.

\begin{proposition} {\rm (p.\ 74, [MUJ]) } Let $U$ be a connected
open subset of a separable Banach space $X$ and let $\{f_n:U \to
\CC \}_{n=1,2,\ldots}$ be a bounded sequence of holomorphic
functions in the compact-open topology.  Suppose also that there
exists a non-empty open subset $V$ of $U$ such that the sequence
$\{f_n (x)\}_n$ converges in $\CC$ for every $x \in V$.  Then, the
sequence $\{f_n\}_n$ converges to a holomorphic function of $U$
uniformly on every compact subset of $U$.
\end{proposition}


Now we turn our attention to characterizations of normal families
that depend on invariant metrics. In what follows, we shall make
use of the Kobayashi metric on a domain $\Omega \subseteq X$.  It
is defined as follows: If $p \in \Omega$ and $\xi \in X$ is a
direction vector then we set
$$
F_K^\Omega(p; \xi) = \inf \left \{ \frac{\|\xi\|}{\|\varphi'(0)\|}
\biggm | \varphi: D \ra \Omega, \varphi(0) = p, \varphi'(0) =
\lambda \xi \ \hbox{for some} \ \lambda \in \RR \right \} \, .
$$
Here $\| \ \ \|$ is the Euclidean norm.

One of the most useful characterizations of normal families, and
one that stems naturally from invariant geometry, is Marty's
criterion.  We now establish such a result in the infinite
dimensional setting.

\begin{proposition} \sl
Let $X$ be a Banach space.  Let $\Omega \ss X$ be a domain and let
${\cal F} = \{f_\alpha\}_{\alpha \in {\cal A}}$ be a family of
holomorphic functions.  The family ${\cal F}$ is normal if and
only if there is a constant $C$ such that, for each (unit)
direction $\xi$,
$$
\frac{|D_\xi f_\alpha(z)|}{1 + |f_\alpha(z)|^2}
  \leq C \cdot F_K^\Omega (z;\xi) \, .
$$
Here $D_\xi$ denotes the directional derivative in the direction
$\xi$.
\end{proposition}
{\bf Proof:}  The proof follows standard lines.  See the proof of
Proposition 1.3 in [CIK, p.\ 306].
\endpf

We next present a rather natural characterization of normal
families that relates the situation on the ambient space to that
on one-dimensional slices (more aptly, one-dimensional analytic
discs):

\begin{proposition} \sl
Let $X$ be a Banach space.  Let $\Omega \ss X$ be a domain and let
${\cal F} = \{f_\alpha\}_{\alpha \in {\cal A}}$ be a family of
holomorphic functions.  The family ${\cal F}$ is normal if and
only if the following condition holds:
\begin{quote}
For each sequence $\varphi_j: D \ra \Omega$ of holomorphic
mappings and each sequence of indices $\alpha_j \in {\cal A}$, $j
= 1,2,\dots$, the family $f_{\alpha_j} \circ \varphi_j$ is normal
on the unit disc $D$. \hfill $(*)$
\end{quote}
\end{proposition}
{\bf Proof:}  The implication ``${\cal F}$ normal $\Rightarrow$
$(*)$'' is immediate from Marty's characterization of normal
families.

For the converse, notice that if Condition $(*)$ holds then, for
each sequence $\varphi_j$ of mappings and each collection
$f_{\alpha_j}$ the compositions $f_{\alpha_j} \circ \varphi_j$
satisfy the conclusion of Marty's theorem:
$$
\frac{[f_{\alpha_j} \circ \varphi_j]'(\zeta)}{1 + |f_{\alpha_j}
\circ \varphi_j|^2}
  \leq C \cdot \frac{1}{1 - |\zeta|^2} \, .   \eqno (\star)
$$
Here the constant $C$ depends in principle on the choice of
$\varphi_j$ and also on the choice of $f_{\alpha_j}$. But in fact
a moment's thought reveals that the choice of $C$ can be taken to
be independent of the choice of these mappings, otherwise there
would be a sequence for which $(\star)$ fails (this is just an
exercise in logic).

But then, using the chain rule, we may conclude that Marty's
Criterion for holomorphic families on a Banach space holds for the
family ${\cal F}$ (see also the proof of Proposition 1.4 in [CIK,
p.\ 307]).  As a result, ${\cal F}$ is normal.
\endpf

%

\section{Other Characterizations of Normality}

It is an old principle of Bloch, enunciated more formally by
Abraham Robinson and actually recorded in mathematical notation by
L. Zalcman (see [ZAL1]), that any ``property'' that would tend to
make an entire function constant would also tend to make a family
of functions normal.  Zalcman's formulation, while incisive, is
rather narrowly bound to the linear structure of Euclidean space.
The paper [ALK] finds a method for formulating these ideas that
will even work on a manifold.  Unfortunately, we must note that
the paper [ALK] has an error, which was kindly pointed out to us
by the authors of [HTT].  We shall include their {\it correct}
formulation of the theorem, and also provide an indication of
their proof.


\begin{proposition} \sl
Let $X$ be a separable Banach space and let $\Omega \ss X$ be a
hyperbolic domain (i.e., a domain on which the Kobayashi metric is
non-degenerate).  Let $Y$ be another separable Banach space.  Let
${\cal F} = \{f_\alpha\}_{\alpha \in {\cal A}} \subseteq
\hbox{Hol}\,(\Omega, Y)$.  The family ${\cal F}$ is {\it not}
normal if and only if there exists a sequence $\{p_j\} \ss \Omega$
with $p_j \ra p_0 \in \Omega$, a sequence $f_j \in {\cal F}$, and
$\{\rho_j\} \ss \RR$ with $\rho_j > 0$ and $\rho_j \ra 0$ such
that
$$
g_j(\xi) = f_j(p_j + \rho_j \xi), \quad \xi \in X
$$
satisfies one of the following assertions:
\begin{enumerate}
\item[{\bf(i)}]  The sequence $\{g_j\}_{j \ge 1}$ is compactly
divergent on $\Omega$;
\item[{\bf (ii)}]  The sequence $\{g_j\}_{j \geq 1}$ converges
uniformly on compact subsets of $\Omega$ to a non-constant
holomorphic mapping $g: \Omega \ra Y$.
\end{enumerate}
\end{proposition}

\begin{remark}  \rm
The error in [ALK] is that the authors did not take into account
the compactly divergent case in the theorem.  Consider the example
(also from [HTT]) of the family ${\cal F}$ of mappings $f_j: D \ra
\CC^2$ given by
$$
f_j(\zeta) = (\alpha_j, \zeta) \, ,
$$
where $1 > \alpha_j > 0$ and $\alpha_j \ra 0$.  Then the family
${\cal F}$ is not normal, but ${\cal F}$ also does not satisfy the
conclusions of part {\bf (ii)} of Proposition 2.1 above, which is
the sole conclusion of the Aladro/Krantz theorem.
\end{remark}

\noindent {\bf Sketch of the Proof of Proposition 2.1:} We first
need a definition. We say that a non-negative, continuous function
$E$ defined on the tangent bundle $TY$ if a {\it length function}
if it satisfies
\begin{enumerate}
\item[{\bf (a)}]  $\displaystyle E(v) = 0 \ \ \hbox{iff} \ \
     v = 0$;
\item[{\bf (b)}]  $\displaystyle E(\alpha v) = |\alpha| E(v) \ \
\hbox{for all} \ \alpha \in \CC \ \hbox{and all} \ v \in TX$.
\end{enumerate}
Now we have:  Let ${\cal F} \ss \hbox{Hol}\,(\Omega, Y)$.  Then
\begin{quote}
\begin{enumerate}
\item[{\bf (1)}]  If ${\cal F}$ is normal then, for each
length function $E$ on $Y$, and for each compact subset $K$ of
$\Omega$, there is a constant $c_K > 0$ such that
$$
E(f(z), df(z)\xi) \leq c_K \cdot \|f\| \quad \hbox{for all} \ z
\in K, \xi \in X \sm \{0\}, f \in {\cal F} \, ; \eqno (*)
$$
\item[{\bf (2)}]  If $Y$ is complete and the family ${\cal F}$ is not
compactly divergent and
satisfies $(*)$ then ${\cal F}$ is normal.
\end{enumerate}
\end{quote}
This result is standard and can be found in [WU] or [HTT]. Now we
treat the result by cases:
\smallskip \\

\noindent {\bf Necessity}
\begin{description}
\item[{\bf Case 1.  The family ${\cal F}$ is compactly divergent.}]
We treat this case in some details since it is new and does not
appear in [ALK]. There is a sequence $\{f_j\} \ss {\cal F}$ that
is compactly divergent. Take $p_0 \in \Omega$ and $r_0 > 0$ such
that $B(p_0, r_0) \subset \! \subset \Omega$. Take $p_j = p_0$ for
all $j \geq 1$ and $\rho_j > 0$ for all $j \geq 1$ such that
$\rho_j \ra 0^+$ and define
$$
g_j(\xi) = f_j(p_j + \rho_j \xi) \, , \quad \hbox{all} \ j \geq 1
\, .
$$
Observe that each $g_j$ is defined on
$$
S_j = \biggl \{\xi \in X: \|\xi\| \leq R_j = \frac{1}{\rho_j}
\hbox{dist} (p_0,
     \partial \Omega) \biggr \} \, .
$$
If $K \ss X$ is compact and $L$ is a compact subset of $Y$ then
there is an index $j_0 \geq 1$ such that $p_0 + \rho_j K \ss
B(p_0, r_0)$ for all $j \geq j_0$.  This implies that $g_j(K) \ss
f_j(\overline{B}(p_0, r_0))$ for each $j \geq j_0$.  Since the
sequence $\{f_j\}$ is compactly divergent, there is an index $j_1
> j_0$ such that $f_j(\overline{B}(p_0, r_0)) \cap L = \emptyset$
for all $j \geq j_1$.  Thus $g_j(K) \cap L = \emptyset$ for all $j
\geq j_1$.  This means that the family $\{g_j\}$ is compactly
divergent.

\item[{\bf Case 2.  The family ${\cal F}$ is not compactly divergent.}]
This follows standard lines, as indicated in [ALK].
\end{description}

\noindent {\bf Sufficiency}
\begin{description}
\item[{\bf Case 1.  The sequence \boldmath $g_j \ra g$ with $g$ not
a constant function.}]  By direct estimation, one shows that
$$
\lim_{j \ra \infty} E(g_j(\xi), dg_j(\xi)(t)) = E(g(\xi),
dg(\xi)(t)) = 0
$$
for $\xi, t \in X$.  Hence $g' \equiv 0$ and so $g$ is constant, a
clear contradiction.  So the family ${\cal F}$ cannot be normal.

\item[{\bf Case 2.  The sequence \boldmath $\{g_j\}$ is compactly
divergent.}]
We may assume that $\{f_j\} \ss {\cal F}$ and $f_j \ra f$.  For
$\xi \in X$ we then have
$$
g_j(\xi) = f_j(p_j + \rho_j \xi) \ra f(p_0) \in Y
$$
since $\rho_j \ra 0$.  This implies that the family $\{g_j\}$ is
{\it not} compactly divergent, a clear contradiction.
\end{description}

That completes our outline of the proof of Proposition 2.1.
\endpf

Constantin Carath\'{e}odory produced a geometric characterization
of normal families that is quite appealing (see [SCH, p.\ 68]). It
has never been adapted even to finitely many complex variables. We
take the opportunity now to offer an infinite dimensional version
(which certainly specializes down to any finite number of
dimensions).

We begin with a little terminology.  Let $\Omega_j$ be domains in
a separable Banach space $X$.  If some Euclidean ball $B(0,r)$, $r
> 0$, is contained in all the domains $\Omega_j$, then
$\hbox{ker}\{\Omega_j\}$ is the largest domain containing 0 and so
that every compact subset of $\hbox{ker}\{\Omega_j\}$ lies in all
but finitely many of the $\Omega_j$.  We say that $\{\Omega_j\}$
{\it converges} to $\Omega_0 \equiv \hbox{ker}\{\Omega_j\}$,
written $\Omega_j \ra \Omega_0$, if every subsequence
$\{\Omega_{j_k}\}$ of these domains has the property that
$\hbox{ker}\{\Omega_{k_j}\} = \Omega_0$.

\begin{theorem} \sl
Fix a separable Banach space $X$.  Let $\{f_n\}$ be a sequence of
univalent, holomorphic mappings from the unit ball ${\bf B}
\subset X$ to another separable Banach space $Y$ with the
properties that
\begin{enumerate}
\item $\displaystyle f_n(0) = 0$;
\item $\displaystyle \langle df_n(0) \11, \11 \rangle > 0$ \, .
\end{enumerate}
[Here $\11$ is the unit vector $(1,0,0,\dots)$.] Set $\Omega_n
\equiv f_n({\cal B})$, $n = 1, 2, \dots$.  Then the $f_n$ converge
normally in ${\cal B}$ to a univalent function $f$ if and only if
\begin{enumerate}
\item $\Omega_0 = \hbox{ker}\{\Omega_n\}$ is hyperbolic and is not
$\{0\}$.
\item $\Omega_n \rightarrow \Omega_0$
\item $\Omega_0 = f({\cal B})$ \, .
\end{enumerate}
\end{theorem}
{\bf Sketch of Proof:}  We first establish that it is impossible
for $\hbox{ker}\{\Omega_n\} = \{0\}$.  Consider the Kobayashi
metric ball ${\bf B} = B_{\cal B}(0,1)$.  Then
$$
f_n: B \rightarrow B_{\Omega_n}(0,1) \, ,
$$
since of course each $f_n$ will be a Kobayashi isometry onto its
image.  Assume that $f_n \ra f$ normally (i.e., uniformly on
compact sets) in ${\cal B}$.  Clearly, under the hypothesis that
$\hbox{ker}\{\Omega_n\} = \{0\}$, there is no $\epsilon > 0$ such
that $b(0,\epsilon) \subseteq \Omega_n$ for $n$ large.  Here $b$
denotes a Euclidean ball.  Thus $B_{\Omega_n}$ must shrink to a
set with no interior. It follows that the sequence $\{f_n\}$
collapses any compact subset of $B$ to a set without interior.
Thus $df \equiv 0$ on $B$ hence $f$ is identically constant. Since
$f(0) = 0$, $f \equiv 0$ (a clear contradiction).

Now we begin proving the theorem proper.  Suppose that $f_n \ra f$
normally on ${\cal B}$ with $f$ univalent. Since, by the preceding
paragraph, $f$ is not identically 0, we may conclude that
$\Omega_0 = \hbox{ker}(\Omega_n) \ne \{0\}$.
\begin{quote}
{\bf CLAIM:}  $\displaystyle f({\cal B}) = \Omega_0.$
\end{quote}
It would follow from this claim that $\Omega_0 \ne X$, for if
$\Omega_0 = X$ then $f^{-1}: X \ra {\cal B}$ univalently,
violating Liouville's theorem [MUJ, p.\ 39].  [It would also
contradict the distance-decreasing property of the Kobayashi
metric.]  Since every subsequence of $\{f_n\}$ converges to $f$,
it follows that every subsequence of $\{\O_n\}$ has kernel $\O_0$.
We write $\O_n \ra \O_0$.
\begin{quote}
{\bf SUBCLAIM I:}  $\displaystyle f({\cal B}) \ss \O_0$.
\end{quote}
For consider any closed metric ball $\overline{B}(0,R) \ss f({\cal
B})$.  We may restrict attention to any finite-dimensional slice
$L$ of this ball, which will of course be compact.  Then $f_n
\biggr |_L \ra f \biggr |_L$. Thus $f_n \biggr |_{L \cap \partial
{\cal B}} \ra f
  \biggr |_{L \cap \partial {\cal B}}$.  As a result,
for $n$ large, we apply the argument principle to any curve in $L
\cap \partial {\cal B}$ to see that each value in
$f(\overline{B}(0,R))$ is attained just once by $f_n$ for $n$
large.  But this just says that $f({\cal B}) \subseteq \O_0$.
\begin{quote}
{\bf SUBCLAIM II:}  $\Omega_0 \subseteq f({\cal B})$.
\end{quote}
For consider $\Omega_0 \ne \{0\}$, and assume $\Omega_0$ is
hyperbolic.  Let $\Omega_n \ra \O_0$.  If $b(0,\epsilon) \ss \O_n$
for all $n$ large, then
$$
\overline{b}(0,\epsilon) \ss B_{\O_n}(0,R) \ss \O_n \quad
\hbox{for} \ n \ \hbox{large}.
$$
So
$$
f_n: B_{\cal B}(0,R) \ra B_{\O_n}(0,R) \supset b(0,\epsilon) \, .
$$
Hence we have a bound from below on the eigenvalues of $d f_n$.

Obversely, we also claim that the eigenvalues of $d f_n$ are
bounded above.  If not, then there exist (Euclidean) unit vectors
$\xi_n$ such that
$$
d f_n (\xi_n) \ra \infty \, .
$$
After a rotation and passing to a subsequence, we can assume that
the $\xi_n$ all point in the direction $\11$.  The result would
then be that $\O_0$ cannot be hyperbolic, a contradiction.

Thus the $\{f_n\}$ are locally bounded and $\{f_n\}$ forms a
normal family, as required.  Thus some subsequence converges (by
the argument principle) to a univalent $f$ such that $f(0) = 0$.
\endpf

\section{A Budget of Counterexamples}

We interrupt our story to provide some examples that exhibit the
limitations of the theory of normal families in infinitely many
variables.

\begin{example} \rm
{\it There is no Montel theorem for holomorphic mappings of
infinitely many variables.}  Indeed, let ${\bf B}$ be the open
unit ball in the Hilbert space $\ell_2$. Define
$$
\varphi_j(\{a_m\}) = \left ( \frac{\sqrt{3/4} a_1}{1 - a_j/2} \ ,
\
  \frac{\sqrt{3/4} a_2}{1 - a_j/2} \ , \  \dots \ , \
  \frac{\sqrt{3/4} a_{j-1}}{1 - a_j/2} \ , \
\frac{a_j - 1/2}{1 - a_j/2} \ , \  \frac{\sqrt{3/4} a_{j+1}}{1 -
a_j/2} \ , \  \dots \right ) \, .
$$
Then each $\varphi_j$ is an automorphism of $B$.

Now fix an index $j$.  Let $K = K_j$ be the compact set $\{
(0,0,\dots, 0, \zeta,0,\dots,0): |\zeta| \leq c\}$, where the
non-zero entry is in the $j^{\rm th}$ position and $1/2 < c < 1$
is a constant. Define the point $p \in K$ to be $p = (0,0, \dots,
0, c, 0, \dots)$, where the non-zero entry is in the $j^{\rm th}$
position. Then
$$
\sup_K \|\varphi_j - \varphi_k \| \geq \| \varphi_j(p) -
\varphi_k(p) \| \geq \left | \frac{c - 1/2}{1 - c/2}  - 0 \right |
= \left | \frac{2c - 1}{2 - c} \right | > 0 \, .
$$

As a result, we see that the sequence $\{\varphi_j\}$ can have no
convergent subsequence.  It also cannot have a compactly divergent
subsequence.
\endpf
\end{example}

\begin{example} \rm
{\it There are no taut domains in infinite dimensional space.}
First we recall H. H. Wu's notion of ``taut''.  Let $N$ be a
complex manifold.  We say that $N$ is {\it taut} if, for every
complex manifold $M$, the family of holomorphic mappings from $M$
to $N$ is normal.  We now demonstrate that there are no such
manifolds in infinite dimensions.

We begin by studying the ball ${\bf B}$ in the Hilbert space
$\ell_2$. We let $N = B$ and $M = D$, the disc in $\CC$ (in fact
it is easy to see that, when testing tautness, it always suffices
to take $M$ to be the unit disc).  Consider the mappings
$$
\varphi_j(\zeta) = \left ( 0 \ , \ 0 \ , \  \dots \ , \  0 \ , \
\frac{1}{3} + \frac{\zeta}{4} \ , \  0 \ , \
  \dots \ , \  0 \right ) \, .
$$
Here the non-zero entry is in the $j^{\rm th}$ position.  Then
$$
\left |  \hbox{image}(\varphi_j) - \hbox{image}(\varphi_k) \right
| \geq \sqrt{ \left ( \frac{1}{12} \right )^2 + \left (
\frac{1}{12} \right )^2} = \sqrt{ \frac{1}{72}} > 0 \, .
$$
Also
$$
\hbox{dist} ( \hbox{image}(\varphi_j), \partial D) = \frac{5}{12}
> 0 \, .
$$
As a result, the sequence $\{\varphi_j\}$ has no convergent
subsequence and no compactly divergent subsequence.

Of course the same argument shows that there is no {\it taut
domain} in Hilbert space, nor is there any taut Hilbert manifold.
\endpf
\end{example}

Of course it should be noted that the Arzela-Ascoli theorem will
fail for families of functions (mappings) taking values in an
infinite dimensional space.  For example, if $X$ is the separable
Hilbert space $\ell_2$ and $f_j: X \ra X$ is given by
$f_j(\{x_j\}) = x_j$ then the $f_j$ are equicontinuous and
equibounded on bounded sets, yet no compact set supports a
uniformly convergent subsequence. Thus the preceding examples do
not come as a great surprise.

It is worth noting that there are results for weak or weak-$*$
normal families that can serve as a good substitute when the
regular (or strong) Montel theorem fails.  We explore some of
these in Section 6.

\section{Normal Functions}

Normal functions were created by Lehto and Virtanen in [LEV] as a
natural context in which to formulate the Lindel\"{o}f principle.
Recall that the Lindel\"{o}f principle says this

\begin{theorem}[Lindel\"{o}f]  \sl
Let $f$ be a bounded holomorphic function on the disc $D$. If $f$
has radial limit $\ell$ at a point $\xi \in \partial D$ then $f$
has non-tangential limit $\ell$ at $\xi$.
\end{theorem}

Lehto and Virtanen realized that boundedness was too strong a
condition, and not the natural one, to guarantee that
Lindel\"{o}f's phenomenon would hold.  They therefore defined the
class of normal functions as follows:

\begin{definition} \rm
Let $f$ be a holomorphic (meromorphic) function on the disc $D \ss
\CC$. Suppose that, for any family $\{\varphi_j\}$ of conformal
self-maps of the disc it holds that $\{f \circ \varphi_m\}$ is a
normal family.  Then we say that $f$ is a normal function.
\end{definition}

Clearly a bounded holomorphic function, a meromorphic function
that omits three values, or a univalent holomorphic function (all
in one complex dimension) will be normal according to this
definition.

Unfortunately, the original definition given by Lehto and Virtanen
is rather limited.  One-connected domains in $\CC^1$ have compact
automorphism groups; finitely connected domains in $\CC^1$, of
connectivity at least two, have finite automorphism group. Generic
domains in $\CC^n$, $n \geq 2$, even those that are topologically
trivial, have automorphism group consisting only of the identity
(such domains are called {\it rigid}). Thus, for most domains in
most dimensions, there are not enough automorphisms to make a
working definition of ``normal function'' possible.  In [CIK],
Cima and Krantz addressed this issue and developed a new
definition of normal function.  We now adapt that definition to
the infinite dimensional case.

\begin{definition} \rm
Let $X$ be a Banach space and let $\Omega$ be a domain in $X$.  A
holomorphic function $f$ on $\Omega$ is said to be {\it normal} if
$$
\frac{|D_\xi f(z)|}{1 + |f(z)|^2} \leq C \cdot F_K^\Omega(z; \xi)
\quad \hbox{for all} \ z \in \Omega, \xi \in X \, .
$$
\end{definition}

\begin{proposition} \sl
Let $f$ be a holomorphic function on a domain $\Omega$ in a Banach
space $X$.  The function $f$ is normal if and only if $f \circ
\varphi$ is normal for each holomorphic $\varphi: D \rightarrow
\Omega$.
\end{proposition}
{\bf Proof:}  The proof is just the same as that in Section 1 of
[CIK].
\endpf

\begin{remark}  \rm
It is a straightforward exercise, using for example Proposition
3.4 (or Marty's characterization of normality), to see that a
holomorphic or meromorphic function on the unit ball ${\bf B}$ in
a separable Hilbert space $H$ is normal if and only if, for every
family $\{\varphi_\alpha\}_{\alpha \in {\cal A}}$ of biholomorphic
self maps of ${\cal B}$, it holds that $\{f \circ
\varphi_\alpha\}$ is a normal family.
\end{remark}

Now let ${\bf B} \ss X$ be the unit ball in a separable Banach
space $X$.  We define a holomorphic function $f$ on ${\cal B}$ to
be {\it Bloch} if
$$
\| df (z) \xi \| \leq C \cdot F_K^\Omega(p; \xi)
$$
for every $z \in {\cal B}$ and every vector $\xi$.  Then it is
routine, following classical arguments, to verify

\begin{proposition} \sl
If $f$ on ${\cal B}$ is a Bloch function then $f$ is normal.
\end{proposition}

\section{Different Topologies on Spaces of Holomorphic Functions}

One way to view a ``normal families'' theorem is that it is a
compactness theorem.  But another productive point of view is to
think of these types of results as relating different topologies
on spaces of holomorphic functions.  We begin our discussion of
this idea by recalling some of the standard topologies, as well as
a few that are more unusual.
\begin{description}
\item[{\bf The Compact-Open Topology}]  In the language
of analysis, this is the topology of uniform convergence on
compact sets.  Certainly in finite-dimensional complex analysis
this is, for many purposes, the most standard topology on general
spaces of holomorphic functions.  In infinite dimensions this
topology is often too coarse (just because compact sets are no
longer very ``fat'').
\item[{\bf The Topology of Pointwise Convergence}]  Here we
say that a sequence $f_j$ of functions or mappings converges if
$f_j(x)$ converges for each $x$ in the common domain $X$ of the
$f_j$.
\item[{\bf The Weak Topology for Distributions}]  Here we
think of a space of holomorphic functions as a subspace of the
space ${\cal E}$ of testing functions for the compactly supported
distributions.  We say that a sequence $f_j$ of holomorphic
functions converges if $\psi(f_j)$ converges for each such
distribution $\psi$.  Of course a similar definition can (and
should) be formulated for nets.
\item[{\bf The Topology of Currents}]  Let $f_j$ be holomorphic
functions and consider the 1-forms $\partial f_j$.  Then we may
think of these forms as {\it currents} lying in the dual of the
space of rectifiable 1-chains; we topologize the $\partial f_j$
accordingly.
\end{description}

\section{A Functional Analysis Approach to Normal Families}

In the classical setting of the unit disc $D \ss \CC$, it is
straightforward to prove that
$$
H^\infty(D) = \left ( L^1(D)/\overline{H^1(D)} \right )^* \, .
\eqno (\star)
$$
Thus $H^\infty$ is a dual space in a natural way.  Properly
viewed, the classical Montel theorem is simply weak-$*$
compactness (i.e., the Banach-Alaoglu theorem) for this dual
space.  Using the Cauchy integral formula as usual, one can see
that convergence in the dual norm certainly dominates uniform
convergence on compact subsets of the disc.

Alternatively, one can think of the elements of $H^\infty(D)$,
with $D$ the disc, as the collection of all operators (by
multiplication) on $H^2$ that commute with multiplication by $z$.
This was Beurling's point of view.  The operator topology turns
out to be equivalent (although this is non-trivial to see) to the
weak-$*$ topology as discussed in the last paragraph, and this in
turn is equivalent to the classical sup-norm topology on
$H^\infty$.

The classical arguments go through to show that there is still a
Beurling theorem on the unit ball in Hilbert space.  It is a
purely formal exercise to verify that $(\star)$ still holds on the
unit ball in $\ell_2$, our usual separable Hilbert space.  As a
result, one can think of the Montel theorem even in infinite
dimensions either in the operator topology or as an application of
the Banach-Alaoglu theorem to $H^\infty$, thought of as a dual
space.

\font\teneufm=eufm10 \font\seveneufm=eufm7 \font\fiveeufm=eufm5
\newfam\eufmfam
\textfont\eufmfam=\teneufm \scriptfont\eufmfam=\seveneufm
\scriptscriptfont\eufmfam=\fiveeufm
\def\frak#1{{\fam\eufmfam\relax#1}}

\newfam\msbfam
\font\tenmsb=msbm10 scaled \magstep1   \textfont\msbfam=\tenmsb
\font\sevenmsb=msbm7 scaled \magstep1 \scriptfont\msbfam=\sevenmsb
\font\fivemsb=msbm5 scaled \magstep1
\scriptscriptfont\msbfam=\fivemsb
\def\Bbb{\fam\msbfam \tenmsb}

\def\RR{{\Bbb R}}
\def\CC{{\Bbb C}}
\def\QQ{{\Bbb Q}}
\def\NN{{\Bbb N}}
\def\ZZ{{\Bbb Z}}
\def\II{{\Bbb I}}
\def\11{{\bf 1}}

\def\e{\epsilon}
\def\q{{\bf q}}
\def\v{{\bf v}}
\def\ra{\rightarrow}
\def\ss{\subseteq}

\def\F{{\cal F}}
\def\O{\Omega}
\def\bO{\partial \Omega}
\def\A{{\cal A}}
\def\endpf{\hfill $\Box$ \smallskip \\}
\def\a{\alpha}
\def\tg{\widetilde{\gamma}}

Now we would like to present an effective weak-normal family
theorem in the context.

Let $Z$ be a Banach space and let $Y = Z*$ be its dual Banach
space.  Let $\Omega$ be an open subset of a Banach space $X$ and
let ${\mathcal O}(\Omega, Y)$ be the set of all holomorphic
mappings from $\Omega$ into $Y$.  Then we consider the topology on
${\mathcal O}(\Omega, Y)$ generated by the sub-basic open sets
given by
$$
G(K,U) \equiv \{f \in {\mathcal O}(\Omega, Y) \mid f(K) \subset
U\}
$$
where $K$ is a compact subset of $\Omega$ and $U$ a weak-* open
subset of $Y$.  Let us call this topology the {\it
compact-weak*-open topology}.

\begin{theorem}  Let $\Omega$ be a domain in a separable Banach
space $X$.  Let $Z$ be a separable Banach space with a countable
Schauder basis, and let $Y = Z^*$.  Further, let $W$ be a bounded
domain in $Y$.  Then the compact-weak*-open topology is a Montel
topology.  In particular, the family ${\mathcal O}(\Omega, Y)$ is
normal with respect to the compact-weak*-open topology.
\end{theorem}

\bf Proof. \rm Let $\{f_j\mid j=1,2,\ldots\} \subset {\mathcal
O}(\Omega, W)$ be given.  We would like to show that there exists
a subsequence that converges in the compact-weak*-open topology.

Let $\{{\bf e}_j \mid j=1,2,\ldots\}$ be a Schauder basis for $Z$.
For $z \in Y$, we define the linear functional $\psi_k : Y \to
\CC$ by $\psi_k (z) = z({\bf e}_k)$. Now we define
$$
\Psi_{k,j} \equiv \psi_k \circ f_j.
$$
Then, we see for each $k$ that the sequence $\{\Psi_{k,j}\}_j$ is
normal by Theorem 1.8.  Therefore we may select subsequences
inductively so that
\begin{itemize}
\item[(1)] $\{\Psi_{1,\sigma_1 (j)} \}_{j=1}^\infty$ is a
subsequence of $\Psi_{1,j}$ which converges in the
compact-weak*-open topology, and
\item[(2)] $\{f_{\sigma_{k+1} (j)}\}_{j=1}^\infty$ is a
subsequence of $\{f_{\sigma_k (j)} \}_{j=1}^\infty$ for every
$k=1,2, \ldots$.
\end{itemize}
Notice that the diagonal sequence $\Psi_{k,\sigma_k(k)} = \psi_k
\circ f_{\sigma_k(k)}$ ($k=1,2,\ldots$) converges in the
compact-weak*-open topology.  Since the weak-* topology separate
points, we may denote the weak-* limit of the sequence
$f_{\sigma_k(k)}(z)$ by $f(z)$ for each $z \in \Omega$. Then the
map $f: \Omega \to Y$ is Gateaux holomorphic.  Since the range of
$f$ is bounded, it follows that $f$ is in fact holomorphic. Thus
the proof is complete. \endpf

\bigskip
Notice that this theorem works for the mappings from the spaces
$\ell^p$ or $c_0$ into the space $\ell^\infty$, for each $p$ with
$1\le p < \infty$.  Therefore, this may be useful for a
characterization problem of infinite dimensional polydisc by its
automorphism group in the space $c_0$ of sequences of complex
numbers converging to zero, for instance.  On the other hand, not
only is this theorem a generalization of the weak-normal family
theorems in the works of Kim/Krantz and Byun/Gaussier/Kim, it also
provides an easier and shorter proof even in the case of separable
Hilbert spaces. See [KIK] and [BGK].
\bigskip

We conclude this section with some examples, due to Jisoo Byun
[BYU], that suggest some of the limitations of normal families in
the infinite dimensional setting.   These examples all relate to
the failure of convexity.

Let $\Omega_1$ and $\Omega_2$ be bounded domains in a Banach space
$X$. We point out that for the holomorphic weak-* limit mapping
$\hat f:\Omega \to X$ of a sequence of holomorphic mappings $f_j :
\Omega_1 \to \Omega_2$ may in general show a surprising behavior
in contrast with the finite dimensional cases.   In the finite
dimensional cases, $\hat f (\Omega_1)$ should be contained in the
closure of $\Omega_2$.  Here we demonstrate that weak-* closure is
about the best one can do, even with the nicest candidates such as
sequences of biholomorphic mappings from the ball.

\begin{example} \rm
  Let $B$ be the unit open ball in $\ell^2$.
Let $\{{\bf e}_j \mid j=1,2, \ldots\}$ be the standard orthonormal
basis for $\ell^2$.  Let $f_k : B \to \ell^2$  be defined by
$$
f_k ({\bf z}) = \sum_{j=1}^{k-1} z_j  {\bf e}_j + (z_k + z_1^2)
{\bf e}_k + \sum_{j=k+1}^\infty z_j {\bf e}_j
$$
where ${\bf z} = z_1 {\bf e}_1 + \ldots$.  Notice that none of
$f_k (B)$ is convex.  In fact, the ball centered at
$\frac{77}{80}{\bf e}_1$ with radius $1/100$ never meets $f_k(B)$,
while it is obvious that the origin and the point ${\bf e}_1$ are
clearly in the norm closure of the union of $f_k(B)$.  Moreover,
the weak limit $\hat f$ of the sequence $f_k$ is the identity map.
Hence $\hat f (B)= B$, which is convex.  This shows that the weak
limit can gain in its image more than the norm closure of the
union of the images of $f_k (B)$.
\end{example}

\begin{example} \rm
  In general the weak-* limit does not make the range convex,
automatically.  If one considers $g_k : B \to \ell^2$ defined by
$$
g_k ({\bf z}) (z_1+z_2^2) {\bf e}_1 + \sum_{j=2}^{k-1} z_j {\bf
e}_j + (z_k + \frac12 z_2^2) {\bf e}_k + \sum_{j=k+1}^\infty z_j
{\bf e}_j,
$$
for $k=3,4,\ldots$.  Then each $g_k$ and the weak limit
$$
\hat g ({\bf z}) = (z_1 + z_2^2) {\bf e}_1 + \sum_{j=2}^\infty z_j
{\bf e}_j
$$
are biholomorphisms of the ball $B$ onto its image.  Notice that
$\hat g(B)$ is not convex.
\end{example}

\section{Many Approaches to Normal Families}

It is natural to try to relate the infinite-dimensional case to
the well-known case of finite dimensions.  In particular, let $\F$
be a family of holomorphic functions on a domain $\O$ in a Banach
space $X$.  Is it correct to say that $\F$ is normal if and only
if the restriction of $\F$ to any finite-dimensional subspace is
normal?  Obversely, if the post-composition of the elements of
$\F$ with each finite-dimensional subspace projection operator is
normal then can we conclude that $\F$ is normal?  We would like to
treat some of these questions here.

\begin{example} \rm
Suppose that if $\F$ is a family of maps of a domain $\Omega$ in a
separable Hilbert space $H$, and assume that
$$
\{\pi_j \circ f: f \in \F\}
$$
is normal for each $\pi_j: H \ra H_j$ the projection of $H$ to the
one-dimensional subspace $H_j$ spanned by the unit vector in the
$j^{\rm th}$ direction.  Then it does not necessarily follow that
$\F$ is a normal family.

To see this, let $H = \ell_2$, and let $f_j(\{x_\ell\}) = x_j$.
Consider each $f_j$ as a map from the unit ball ${\bf B} \subseteq
H$ to itself.  Then, for each fixed $k$, $\{\pi_k \circ
f)j\}_{j=1}^\infty$ is a normal family, yet the family $\F =
\{f_j\}_{j=1}^\infty$ is definitely not normal.
\end{example}

\begin{example} \rm
Suppose that if $\F$ is a family of maps of a domain $\Omega$ in a
separable Hilbert space $H$. Suppose that, for each $k$, the
collection
$$
\{f \circ \mu_k: f \in \F\}
$$
is normal for each $\mu_j: \CC \ra H$ the injection of $\CC$ to
$H$ in the $j^{\rm th}$ variable.  Then it does not necessarily
follow that $\F$ is a normal family.

To see this, again consider $H = \ell_2$, and let $f_j(\{x_\ell\})
= x_j$. Consider each $f_j$ as a map from the unit ball ${\bf B}
\subseteq H$ to itself.  Then, for each fixed $k$, the family $\{f
\circ \mu_k\}_{f \in \F}$ is clearly normal.  Yet the entire
family $\F$ is plainly not normal---as we discussed in Example
2.1.  The reader should compare this example to Proposition 1.14,
which gives a positive result along these lines.
\end{example}

One of the main lessons of the classic paper [WU] by H. H. Wu is
that the normality or non-normality of a family of mappings
depends essentially on the target space (this is the provenance of
the notion of {\it taut} manifold).  With this point in mind, we
now formulate a counterpoint to Example 2.1:

\begin{proposition}  \sl
Let ${\cal C} = \{ \{x_j\}_{j=1}^\infty : |x_j| \leq 1/j\}$ be the
Hilbert cube.  Let $H = \ell_2$ be the canonical separable Hilbert
space.  Then any family $\F$ from a domain $\Omega \ss H$ to
${\cal C}$ will be normal.
\end{proposition}
\noindent {\bf Proof:}  It suffices to prove that the correct
formulation of the Arzela-Ascoli theorem holds.  In particular, we
establish this result:
\begin{quote}
If ${\cal G} = \{g_\alpha\}_{\alpha \in {\cal A}}$ is a family of
functions from a domain $\Omega \ss H$ into ${\cal C}$ which is
{\bf (i)}  equibounded and {\bf (ii)} equicontinuous, then ${\cal
G}$ has a uniformly convergent subsequence.
\end{quote}
In fact the usual proof of Arzela-Ascoli, that can be found in any
text (see, for instance, [KRA, p.\ 284]]), goes through once we
establish this basic fact:  If $g_\alpha: \Omega \ra {\cal C}$ and
$x_0 \in \Omega$ is fixed then $\{g_\alpha(x_0)\}$ has a
convergent subsequence.  Of course this simple assertion is the
consequence of a standard diagonalization argument.
\endpf

\newpage

\noindent {\Large \bf References}
\medskip  \\

\begin{quote}
\begin{enumerate}
\item[{\bf [ALK]}]  G. Aladro and S. G. Krantz, A criterion
for normality in $\CC^n$, {\it Jour.\ Math.\ Anal.\ and Appl.}
161(1991), 1-8.

\item[{\bf [ALM]}]  F. J. Almgren, Mass continuous cochains are
differential forms. {\it Proc.\ Amer.\ Math.\ Soc.} 16(1965),
1291--1294.

\item[{\bf [ACP]}]  J. M. Anderson, J. Clunie, and Ch. Pommerenke, On Bloch
functions and normal functions, {\it J. Reine Angew.\ Math.}
270(1974), 12--37.

\item[{\bf [ANT]}]  J. Ansemil and J. Taskinen, On a problem of
topologies in infinite-dimensional holomorphy. {\it Arch. Math.}
(Basel) 54 (1990), no. 1, 61--64.

\item[{\bf [ARA]}] J. Arazy, An application of
infinite-dimensional holomorphy to the geometry of Banach spaces.
Geometrical aspects of functional analysis (1985/86), 122--150,
{\it Lecture Notes in Math.}, 1267, Springer, Berlin-New York,
1987.

\item[{\bf [BAR]}]  T. J. Barth, Separate analyticity, separate normality,
and radial normality for mappings, {\it Several Complex Variables
(Proc.\ Sumpos.\ Pure Math., Vol.\ XXX, Part 2, Williams College,
Williamstown, Mass., 1975)}, pp.\ 221--224.  American Mathematical
Society, Providence, R.I., 1977.

\item[{\bf [BAR]}]  T. J. Barth, Normality domains for families of
holomorphic maps, {\it Math.\ Annalen} 190 (1971), 293--297.

\item[{\bf [BAR]}]  T. J. Barth, Extension of normal families of
holomorphic functions, {\it Proc.\ Amer.\ Math.\ Soc.} 16(1965),
1236--1238.

\item[{\bf [BJL]}]  S. Bjon and M. Lindstr\"{o}m, On a bornological
structure in infinite-dimensional holomorphy. {\it Math.\ Nachr.}
139(1988), 77--86.

\item[{\bf [BMN1]}]  J. Barroso, M. Matos, and L. Nachbin,
On holomorphy versus linearity in classifying locally convex
spaces. {\it Infinite dimensional holomorphy and applications}
(Proc. Internat. Sympos., Univ. Estadual de Campinas, S ao Paulo,
1975), pp. 31--74. North-Holland Math.  Studies, Vol. 12, Notas de
Mat., No. 54, North-Holland, Amsterdam, 1977.

\item[{\bf [BMN2]}]  J. Barroso, M. Matos, and L. Nachbin,
On bounded sets of holomorphic mappings. Proceedings on Infinite
Dimensional Holomorphy (Internat. Conf., Univ. Kentucky,
Lexington, Ky., 1973), pp. 123--134. {\it Lecture Notes in Math.},
Vol. 364, Springer, Berlin, 1974.

\item[{\bf [BAN1]}]  J. Barroso and L. Nachbin, A direct
sum is holomorphically bornological with the topology induced by a
Cartesian product.  {\it Portugal. Math.} 40(1981), no. 2,
252--256.

\item[{\bf [BAN2]}]  J. Barroso, Jorge Alberto and L. Nachbin,
Sur certaines proprietes bornologiques des espaces d'applications
holomorphes. (French) {\it Troisieme Colloque sur l'Analyse
Fonctionnelle} (Liege, 1970), pp. 47--55. Vander, Louvain, 1971.

\item[{\bf [BAY]}]  A. Bayoumi, Infinite-dimensional holomorphy without
convexity condition. I. The Levi problem in nonlocally convex
spaces. New frontiers in algebras, groups and geometries
(Monteroduni, 1995), 287--306, {\it Ser.\ New Front.\ Adv.\ Math.\
Ist.\ Ric.\ Base}, Hadronic Press, Palm Harbor, FL, 1996.

\item[{\bf [BJL]}]  S. Bjon and M. Lindstr\"{o}m, A general approach to
infinite-dimensional holomorphy. {\it Monatsh.\ Math.} 101 (1986),
no. 1, 11--26.

\item[{\bf [BOL]}]  P. Boland, An example of a nuclear space in
infinite dimensional holomorphy. {\it Ark.\ Mat.} 15 (1977), no.
1, 87--91.

\item[{\bf [BMV]}]  M. B\"{o}rgens, R. Meise, and D. Vogt,
$\Lambda (\alpha)$-nuclearity in infinite-dimensional holomorphy.
{\it Math.\ Nachr.} 106 (1982), 129--146.

\item[{\bf [BRJ]}] H. Braunss and H. Junek, On types of polynomials and
holomorphic functions on Banach spaces. {\it Note Mat.} 10 (1990),
no. 1, 47--58.

\item[{\bf [BYU]}]  J. Byun, Geometry of automorphism group orbits
and Levi Geometry, Ph.\ D. Thesis, Pohang University of Science
and Technology, Pohang 790-784 Korea, June 2002.

\item[{\bf [BGK]}]  J. Byun, H. Gaussier, and K. T. Kim,
Weak-type normal families of holomorphic mappings in Banach spaces
and characterization of the Hilbert ball by its automorphism
group, {\it J. Geom.\ Analysis}, to appear.

\item[{\bf [CHU]}]  C.-T. Chung, {\it Normal Families of
Meromorphic Functions}, World Scientific, Singapore, 1993.

\item[{\bf [COL1]}]  J. Colombeau, Some aspects of
infinite-dimensional holomorphy in mathematical physics.  {\it
Aspects of mathematics and its applications}, 253--263,
North-Holland Math. Library, 34, North-Holland, Amsterdam-New
York, 1986.

\item[{\bf [COL2]}]  J. Colombeau, On some various notions of infinite
dimensional holomorphy. {\it Proceedings on Infinite Dimensional
Holomorphy} (Internat.  Conf., Univ. Kentucky, Lexington, Ky.,
1973), pp. 145--149. Lecture Notes in Math., Vol. 364, Springer,
Berlin, 1974.

\item[{\bf [COP]}]  J. Colombeau and B. Perrot, Reflexivity and kernels
in infinite-dimensional holomorphy. {\it Portugal.\ Math.} 36
(1977), no. 3-4, 291--300 (1980).

\item[{\bf [CIK]}]  J. A. Cima and S. G. Krantz, The Lindel\"{o}f principle
and normal functions of several complex variables, {\it Duke
Math.\ Jour.} 50(983), 303--328.

\item[{\bf [DIN1]}]  S. Dineen, Monomial expansions in
infinite-dimensional holomorphy. {\it Advances in the theory of
Frechet spaces} (Istanbul, 1988), 155--171, NATO Adv. Sci. Inst.
Ser. C: Math. Phys. Sci., 287, Kluwer Acad. Publ., Dordrecht,
1989.

\item[{\bf [DIN2]}]  S. Dineen, Surjective limits of locally convex
spaces and their application to infinite dimensional holomorphy.
{\it Bull. Soc. Math. France} 103 (1975), no. 4, 441--509.

\item[{\bf [DIN3]}]  S. Dineen, {\it The Schwarz Lemma}, Oxford University
Press, New York, 1989.

\item[{\bf [DIN4]}]  S. Dineen, {\it Complex Analysis of Infinite
Dimensional Spaces}, Springer-Verlag, New York, 1999.

\item[{\bf [DIN4]}]  S. Dineen,  {\it Complex Analysis in Locally
Convex Spaces}, North-Holland, Amsterdam, 1981.

\item[{\bf [GMR]}]  D. Garcia, M. Maestre, P. Rueda,  Weighted spaces
of holomorphic functions on Banach spaces. {\it Studia Math.} 138
(2000), no. 1, 1--24.

\item[{\bf [GRK]}]  L. Gruman and C. Kiselman, Le probleme de Levi
dans les espaces de Banach \`{a} base. (French) {\it C. R. Acad.
Sci. Paris} Ser. A-B 274(1972), A1296--A1299.

\item[{\bf [HAH1]}]  K. T. Hahn, Equivalence of the classical theorems
of Schottky, Landau, Picard and hyperbolicity, {\it Proc.\ Amer.\
Math.\ Soc.} 89(1983), 628--632.

\item[{\bf [HAH2]}]  K. T. Hahn, Higher dimensional generalizations of
the Bloch constant and their lower bounds, {\it Trans.\ Amer.\
Math.\ Soc.} 179(1973), 263--274.

\item[{\bf [HTT]}]  P. D. Huong, P. N. T. Trang, and D. D. Thai,
Families of normal maps in several complex variables and
hyperbolicity of complex spaces, preprint.

\item[{\bf [JAP]}]  Jarnicki and P. Pflug, {\it Invariant Distances and
Metrics in Complex Analysis}, de Gruyter, New York, 1993.

\item[{\bf [KIK]}]  K. T. Kim and S. G. Krantz, Characterization of
the Hilbert ball by its automorphism group, {\it Transactions of
the AMS} 354(2002), 2797--2828.

\item[{\bf [KOB1]}]  S. Kobayashi, {\it Hyperbolic Manifolds and
Holomorphic Mappings}, Dekker, New York, 1970.

\item[{\bf [KOB2]}]  S. Kobayashi, {\it Hyperbolic Complex Spaces},
Springer, New York, 1998.

\item[{\bf [KOB3]}]  S. Kobayashi, Intrinsic distances, measures,
and geometric function theory, {\it Bull.\ Amer.\ Math.\ Soc.}
82(1976), 357--416.

\item[{\bf [KRA]}]  S. G. Krantz, {\it Real Analysis and Foundations},
CRC Press, Boca Raton, Florida, 1991.

\item[{\bf [LEV]}]  O. Lehto and K. Virtanen, Boundary behavior and normal
meromorphic functions, {\it Acta Math.} 97(1957), 47--65.

\item[{\bf [LEM1]}]  L. Lempert, Approximation of holomorphic functions of
infinitely many variables. II. {\it Ann.\ Inst.\ Fourier}
(Grenoble) 50 (2000), no. 2, 423--442.

\item[{\bf [LEM2]}]  L. Lempert,  Approximation de fonctions
holomorphes d'un nombre infini de variables. (French)
[Approximation of holomorphic functions of infinitely many
variables] {\it Ann. Inst. Fourier} (Grenoble) 49 (1999), no. 4,
1293--1304.

\item[{\bf [LEM3]}]  L. Lempert,  The Dolbeault complex in infinite
dimensions. I. {\it J. Amer. Math. Soc.} 11(1998), no. 3,
485--520.

\item[{\bf [LEM4]}]  L. Lempert,  The Dolbeault complex in infinite
dimensions. II. {\it Jour.\ Amer.\ Math.\ Soc.} 12(1999),
775--793.

\item[{\bf [LEM5]}]  L. Lempert,  The Dolbeault complex in infinite
dimensions. III. {\it J.\ Amer.\ Math.\ Soc.} 12(1999), no. 3,
775--793.

\item[{\bf [LEM6]}]  L. Lempert,  The Cauchy-Riemann equations in
infinite dimensions. {\it Journees Equations aux Derivees
Partielles} (Saint-Jean-de-Monts, 1998), Exp. No. VIII, 8 pp.,
Univ. Nantes, Nantes, 1998.

\item[{\bf [LEP]}]  A. Lepskiui, A superholomorphy criterion and
embedding theorems in infinite- dimensional holomorphy. (Russian)
{\it Izv. Vyssh. Uchebn. Zaved. Mat.} 1992, no. 11, 72--75 (1993)
translation in Russian Math. (Iz. VUZ) 36 (1992), no. 11, 70--73.

\item[{\bf [LIR]}]  M. Lindstr\"{o}m and R. Ryan, Applications of
ultraproducts to infinite-dimensional holomorphy. {\it Math.
Scand.} 71(1992), no. 2, 229--242.

\item[{\bf [MAN1]}] M. Matos and L. Nachbin, Reinhardt domains of
holomorphy in Banach spaces. {\it Adv. Math.} 92(1992), no. 2,
266--278.

\item[{\bf [MAN2]}] M. Matos and L. Nachbin, Entire functions
on locally convex spaces and convolution operators. {\it
Compositio Math.} 44(1981), no.  1-3, 145--181.

\item[{\bf [MIN]}]  K. Min and L. Nel, Infinite-dimensional
holomorphy via categorical differential calculus. {\it Monatsh.
Math.} 111 (1991), no. 1, 55--68.

\item[{\bf [MON1]}]  P. Montel, Sur les familles quasi normals
de fonctions holomorphes, {\it Memoires de la Classe des
sciences}, Bruxelles, 1922.

\item[{\bf [MON2]}]  P. Montel, {\it Lecons sur les Familles Normales
de Fonctions Analytiques et Leur Applications}, Gauthier-Villars,
Paris, 1927.

\item[{\bf [MOR]}] L. Moraes, Quotients of spaces of holomorphic
functions on Banach spaces. {\it Proc.\ Roy.\ Irish Acad.} Sect. A
87 (1987), no. 2, 181--186.

\item[{\bf [MUJ]}]   J. Mujica, {\it Complex Analysis in Banach Spaces},
North-Holland, Amsterdam and New York, 1986.

\item[{\bf [MNA]}]  J. Mujica and L. Nachbin, Linearization of
holomorphic mappings on locally convex spaces. {\it J. Math. Pures
Appl.} (9) 71(1992), no.  6, 543--560.

\item[{\bf [NAC1]}]  L. Nachbin, When does finite holomorphy imply
holomorphy? {\it Portugal. Math.} 51 (1994), no. 4, 525--528.

\item[{\bf [NAC2]}]  L. Nachbin, Some aspects and problems in holomorphy.
Extracta Math. 1 (1986), no. 2, 57--72.

\item[{\bf [NAC3]}]  L. Nachbin, On pure uniform holomorphy in spaces of
holomorphic germs. {\it Resultate Math.} 8(1985), no. 2, 117--122.

\item[{\bf [NAC4]}]  L. Nachbin, Why holomorphy in infinite dimensions?
{\it Enseign. Math.} (2) 26 (1980), no. 3-4, 257--269 (1981).

\item[{\bf [NAC5]}]  L. Nachbin, Analogies entre l'holomorphie et la
linearite. (French) {\it Seminaire Paul Kree}, $4^{rm e}$ annee:
1977--1978. equations aux derivees partielles en dimension
infinie, Exp. No. 1, 10 pp., Secretariat Math., Paris, 1979.

\item[{\bf [NAC6]}]  L. Nachbin, Warum unendlichdimensionale
Holomorphie? (German) {\it Jahrbuch \"{U}berblicke Mathematik},
1979, pp. 9--20, Bibliographisches Inst., Mannheim, 1979. x

\item[{\bf [NAC7]}]  L. Nachbin, Some holomorphically significant
properties of locally convex spaces. {\it Functional analysis}
(Proc. Brazilian Math. Soc.  Sympos., Inst. Mat. Univ. Estad.
Campinas, Sao Paulo, 1974), pp. 251--277. Lecture Notes in Pure
and Appl. Math., 18. Dekker, New York, 1976.

\item[{\bf [NAC8]}]  L. Nachbin, Some problems in the application
of functional analysis to holomorphy. Advances in holomorphy ({\it
Proc. Sem.\ Univ.\ Fed.\ Rio de Janeiro}, Rio de Janeiro, 1977),
pp. 577--583. North-Holland Math. Studies, 34. North-Holland,
Amsterdam, 1979.

\item[{\bf [NAC9]}]  L. Nachbin, Limites et perturbation des
applications holomorphes. (French) Fonctions analytiques de
plusieurs variables et analyse complexe ({\it Colloq. Internat.
CNRS}, No. 208, Paris, 1972), pp. 141--158. ``Agora Mathematica'',
No. 1, Gauthier-Villars, Paris, 1974.

\item[{\bf [NAC10]}]  L. Nachbin, Sur quelques aspects recents de
l'holomorphie en dimension infinie. (French) {\it Seminaire
Goulaouic-Schwartz} (1971--1972), equations aux derivees
partielles et analyse fonctionnelle, Exp. No. 18, 9 pp. Ecole
Polytech., Centre de Math., Paris, 1972.

\item[{\bf [NAC11]}]  L. Nachbin, A glimpse at infinite dimensional
holomorphy. Proceedings on Infinite Dimensional Holomorphy,
(Internat. Conf., Univ.  Kentucky, Lexington, Ky., 1973), pp.
69--79. {\it Lecture Notes in Math.}, Vol. 364, Springer, Berlin,
1974.

\item[{\bf [NAC12]}]  L. Nachbin, Concerning holomorphy types for Banach
spaces. Memerias de Matematica da Universidade Federal do Rio de
Janeiro, No.  3. [Mathematical Reports of the Federal University
of Rio de Janeiro, No. 3] Instituto de Matematica, Universidade
Federal do Rio de Janeiro; Coordena ao dos Programas de
Pos-Graduat ao em Engenharia, Universidade Federal do Rio de
Janeiro, Rio de Janeiro, 1971. i+11 pp.

\item[{\bf [NAC13]}]  L. Nachbin, On vector-valued versus scalar-valued
holomorphic continuation. {\it Nederl.\ Akad.\ Wetensch.\ Proc.}
Ser. A $=$ {\it Indag. Math.} 35(1973), 352--354.

\item[{\bf [NAC14]}]  L. Nachbin, Recent developments in infinite
dimensional holomorphy. {\it Bull.\ Amer.\ Math.\ Soc.} 79(1973),
625--640.

\item[{\bf [NAC15]}]  L. Nachbin, Concerning holomorphy types for Banach
spaces. {\it Studia Math.} 38 1970 407--412.

\item[{\bf [NAC16]}]  L. Nachbin, A glimpse at infinite dimensional
holomorphy. {\it Proceedings on Infinite Dimensional Holomorphy}
(Internat. Conf., Univ.  Kentucky, Lexington, Ky., 1973), pp.
69--79. Lecture Notes in Math., Vol. 364, Springer, Berlin, 1974.

\item[{\bf [NAC17]}]  L. Nachbin, Recent developments in infinite
dimensional holomorphy.  {\it Bull.\ Amer.\ Math.\ Soc.} 79
(1973), 625--640.

\item[{\bf [NAC18]}]  L. Nachbin, {\it Topology on Spaces of
Holomorphic Functions}, Springer-Verlag, New York, 1969.

\item[{\bf [RIC]}] C. E. Rickart, A function algebra approach to infinite
dimensional holomorphy. {\it Analyse fonctionnelle et
applications} (Comptes Rendus Colloq.  Analyse, Inst. Mat., Univ.
Fed. Rio de Janeiro, Rio de Janeiro, 1972), pp. 245--260.
Actualites Sci. Indust., No. 1367, Hermann, Paris, 1975.

\item[{\bf [SCH]}]  J. Schiff, {\it Normal Families}, Springer, New York,
1993.

\item[{\bf [UPM]}] H. Upmeier, Some applications of
infinite-dimensional holomorphy to mathematical physics. {\it
Aspects of mathematics and its applications}, 817--832,
North-Holland Math. Library, 34, North-Holland, Amsterdam-New
York, 1986.

\item[{\bf [WAE]}]  L. Waelbroeck, The holomorphic functional
calculus and infinite dimensional holomorphy. {\it Proceedings on
Infinite Dimensional Holomorphy} (Internat. Conf., Univ. Kentucky,
Lexington, Ky., 1973), pp. 101--108. Lecture Notes in Math., Vol.
364, Springer, Berlin, 1974.

\item[{\bf [WU]}]  H. Wu, Normal families of holomorphic mappings,
{\it Acta Mathematica} 119(1967), 193--223.

\item[{\bf [ZAL1]}]  L. Zalcman, A heuristic principle in complex
function theory, {\it Am.\ Math.\ Monthly} 82(1975), 813--817.

\item[{\bf [ZAL2]}]  L. Zalcman, Normal families:  new perspectives,
{\it Bull.\ Amer.\ Math.\ Soc.} 35(1998), 215--230.

\item[{\bf [ZAL3]}]   L. Zalcman, New light on normal families,
{\it Proceedings of the Ashkelon Workshop on Complex Function
Theory (1996)}, 237--245, {\it Israel Math. Conf. Proc.} 11,
Bar-Ilan Univ., Ramat Gan, 1997.

\item[{\bf [ZALD]}] I. Zalduendo, Duality and extensions in
infinite-dimensional holomorphy. (Spanish) {\it Proceedings of the
Second Latin American Colloquium on Analysis} (Spanish) (Santafe
de Bogota, 1992). Rev. Colombiana Mat. 27 (1993), no. 1-2,
131--135.
\end{enumerate}
\end{quote}

\newpage
Kang-Tae Kim

Department of Mathematics

Pohang University of Science and Technology

Pohang 790-784 Korea

{\tt kimkt@postech.ac.kr}
\bigskip

Steven G. Krantz

Department of Mathematics

Washington University

St. Louis, Missouri 63130 U.S.A.

{\tt sk@math.wustl.edu}

\end{document}